# Getting rid of the negative Schwarzian derivative condition

By O. S. Kozlovski

## Abstract

In this paper we will show that the assumption on the negative Schwarzian derivative is redundant in the case of $C^3$ unimodal maps with a nonflat critical point. The following theorem will be proved: For any $C^3$ unimodal map of an interval with a nonflat critical point there exists an interval around the critical value such that the first entry map to this interval has negative Schwarzian derivative. Another theorem proved in the paper provides useful cross-ratio estimates. Thus, all theorems proved only for unimodal maps with negative Schwarzian derivative can be easily generalized.

## 1. Introduction

If a map $f$ has critical points, one cannot hope to get a bound for its nonlinearity. However, if the map has one extra property, namely the negative Schwarzian derivative, then the behavior of this map is somewhat similar to the behavior of univalent maps. For such maps there are analogies to the Koebe lemma or to the minimum modulus principle, but their main property is that they increase cross-ratios. This property appears to be crucial for the whole theory of one-dimensional maps. There are still many theorems which were proved only for maps having negative Schwarzian derivative. To generalize a theorem for maps without negative Schwarzian derivative one would have to estimate lengths of intervals of some orbits. It was not always easy to make these estimates.

Moreover, an assumption on negative Schwarzian derivative is unnatural. Indeed, this negative Schwarzian derivative condition does not have (and cannot have) any dynamical meaning. A smooth change of the coordinate can destroy this property of a map. It was also not clear if the class of $S$-maps (i.e. maps with negative Schwarzian derivative) has some special properties in "small scales" which smooth maps do not enjoy.

The main purpose of this paper is to provide a universal tool which enables one to deduce any statement proved for unimodal maps to the case of smooth maps. So, there is nothing special about $S$-unimodal maps!

On the other hand, maps with negative Schwarzian derivative have many special properties which do not hold for other maps. For example, an $S$-unimodal map can have at most two attracting periodic points (the basin of attraction of one critical point should contain a boundary point of the interval and the basin of attraction of the other one should contain the critical point) while an arbitrary unimodal map can have arbitrarily many attracting periodic points. However, all these extra properties have a global nature and they cannot occur in "small scales".

It is shown in [3dMvS] that all nice properties of the $S$-maps are consequences of one: $S$-maps increase cross-ratios. However, even if the map does not increase cross-ratios but remains bounded from zero, then one can prove analogies to the theorems for the negative Schwarzian derivative (for example, the Koebe principle). And it appears that this bound for the distortion of the cross-ratio exists provided the last interval from the orbit is small. Thus, the next theorems (see Section 3) allow us to transfer the properties of $S$-unimodal maps to the "small scales" of arbitrary $C^3$ unimodal maps. For example, the theorem we mentioned above can be transferred in the following statement: for any $C^3$ unimodal map $f$ with a nonflat critical point the periods of sinks are uniformly bounded.

First, the negative Schwarzian derivative condition in the content of maps of interval was introduced by Singer who noticed that if a map has negative Schwarzian, then all of its iterates have negative Schwarzian as well and that if the Schwarzian derivative of a map $f$ is negative, then $|Df|$ cannot have a positive local minimum, [Sin]. In fact, the Schwarzian derivative has already been used before by Herman, [Her]. Guckenheimer, Misiurewicz and van Strien showed the importance of the Schwarzian derivative for the study of several dynamical properties, [Guc], [Mis], [1vS]. Later, a large number of papers appeared where extensive use of the negative Schwarzian derivative condition was made. For a comprehensive list of these papers see [3dMvS]. Then it was realized that the maps with negative Schwarzian derivative increase some cross-ratios and that it is a very powerful tool; see [Pre], [Yoc], [1dMv], [2dMv] and afterwards [Swi] (see also Lemma 2.3 below). In [2vS] the cross-ratio **a** (see the next section for its definition) was used to analyze situations where one has only detailed information on one side of the orbit of some interval.

I am deeply grateful to S. van Strien for many useful remarks and comments and for his constant interest in my work. I would like to thank J. Graczyk, K. Khanin, G. Levin, M. St. Pierre, G. Świątek, M. Tsujii, and, particularly, W. de Melo for interesting discussions and remarks. This work has been supported by the Netherlands Organization for Scientific Research (NWO).

## 2. Schwarzian derivative and cross ratios.

Before giving statements of the main theorems of the paper we have to define the Schwarzian derivative and cross-ratios to be used.

Let $f$ be a $C^3$ map of an interval. The Schwarzian derivative $Sf$ of the map $f$ is defined for noncritical points of $f$ by the formula:

$$Sf(x) = \frac{D^3 f(x)}{Df(x)} - \frac{3}{2}\left(\frac{D^2 f(x)}{Df(x)}\right)^2.$$

One can easily check the following expression of the Schwarzian derivative of a composition of two maps:

$$S(fg)(x) = Sf(g(x))\,(Dg(x))^2 + Sg(x).$$

From this formula we can deduce an important property of maps having negative Schwarzian derivative (i.e. $Sf(x) < 0$ where $Df(x) \neq 0$): all iterates of such maps also have negative Schwarzian derivative.

Since we cannot control the distortion of a map with critical points, the ratio of lengths of two adjacent intervals can change dramatically under iterates of the map. So, instead of considering three consecutive points, we consider four points and we measure their positions by their cross-ratios. There are several types of cross-ratios which work more or less in the same way. We will use just a standard cross-ratio which is given by the formula:

$$\mathbf{b}(M, J) = \frac{|J||M|}{|M^-||M^+|}$$

where $J \subset M$ are intervals and $M^-$, $M^+$ are connected components of $M \setminus J$.

Another useful cross-ratio (which is in some sense degenerate) is the following:

$$\mathbf{a}(M, J) = \frac{|J||M|}{|M^- \cup J||J \cup M^+|}$$

where the intervals $M^-$ and $M^+$ are defined as before.

If $f$ is a map of an interval, we will measure how this map distorts the cross-ratios and introduce the following notation:

$$\mathbf{B}(f, M, J) = \frac{\mathbf{b}(f(M), f(J))}{\mathbf{b}(M, J)},$$

$$\mathbf{A}(f, M, J) = \frac{\mathbf{a}(f(M), f(J))}{\mathbf{a}(M, J)}.$$

The main property of maps with negative Schwarzian derivative in given in the following well-known theorem:

LEMMA 2.1. *Let $f$ be a $C^3$ map with negative Schwarzian derivative and $M$ be an interval such that $f|_M$ is a diffeomorphism. Then for any subinterval $J \subset M$,*

$$\mathbf{A}(f, M, J) \geq 1,$$
$$\mathbf{B}(f, M, J) \geq 1.$$

In fact, almost all other properties of maps with negative Schwarzian derivative are consequences of this theorem. We will need only the most powerful tool, the so-called Koebe principle, which controls the distortion of maps away from the critical points. If the map satisfies the negative Schwarzian derivative condition, then we can apply the next result to each iterate. Otherwise, the main theorems of this paper allow us to apply it anyway, provided the image of the interval is small enough.

LEMMA 2.2 (the Koebe Principle). *Let $J \subset M$ be intervals, $f : M \to \mathbf{R}$ be a $C^1$ diffeomorphism, $C$ be a constant such that $0 < C < 1$. Assume that for any interval $J^*$ and $M^*$ with $J^* \subset M^* \subset M$,*

$$\mathbf{B}(f, M^*, J^*) \geq C.$$

*If $f(M)$ contains a $\tau$-scaled neighborhood of $f(J)$, then*

$$\frac{1}{K(C, \tau)} \leq \frac{Df(x)}{Df(y)} \leq K(C, \tau)$$

*where $x, y \in J$ and $K(C, \tau) = \frac{(1+\tau)^2}{C^6 \tau^2}$.*

Here we say that an interval $M$ is a $\tau$-scaled neighborhood of the interval $J$, if $M$ contains $J$ and if each component of $M \setminus J$ has at least length $\tau|J|$.

The proofs of Lemmas 2.1 and 2.2 can be found in [3dMvS].

So, one has good nonlinearity estimates if bounds on the distortion of the cross-ratio are known. In this section we will formulate a lemma which describes the distortion of the cross-ratios under high iterates of a smooth map provided some summability conditions are satisfied.

The maps which the next lemma can be applied to should have *a nonflat* critical point. If the map is smooth and one of its higher derivatives does not vanish at the critical point, this map automatically has a nonflat critical point. If the map $f$ is only $C^3$, then we will say that $f$ has a nonflat critical point if there is a local $C^3$ diffeomorphism $\phi$ with $\phi(c) = 0$ such that $f(x) = \pm|\phi(x)|^\alpha + f(c)$ for some real $\alpha \geq 2$. Thus we assume that the order of the critical point is the same on both sides.

Maps which do have a flat critical point, may have completely different properties. For example, such maps can have wandering intervals.

The following result is well-known and can be found in [1dMv], [2dMv] and [2vS]. It is based on the very simple idea: near the nonflat critical point the map has negative Schwarzian derivative and outside of a fixed neighborhood of the critical point the distortion of the map is bounded.

LEMMA 2.3 ([2dMv]). *Let $X$ be an interval, $f : X \to X$ be a $C^{2+1}$ map whose critical points are non-flat. Then there exists a constant $C_1$ with the following property. If $M \supset J$ are intervals such that $f^m$ is a diffeomorphism on $M$ and $M \setminus J$ consists of two components $M^-$ and $M^+$ then:*

$$\mathbf{A}(f^m, M, J) \geq \exp\left\{C_1 \sum_{i=0}^{m-1} |f^i(M^-)| |f^i(M^+)|\right\},$$

$$\mathbf{B}(f^m, M, J) \geq \exp\left\{C_1 \sum_{i=0}^{m-1} |f^i(M)|^2\right\}.$$

The negative Schwarzian derivative condition can be introduced for $C^2$ maps as well, but we will not consider such maps here and all results of this paper about the negative Schwarzian derivative condition can be applied to the case of $C^{2+1}$ maps.

## 3. How to get rid of the negative Schwarzian derivative condition

Here we give the main theorems of the paper which we formulate and prove only in the unimodal case, i.e. for maps of an interval which have only one turning point. However, it seems that the method can be applied to the multimodal case as well; to prove analogous theorems in the multimodal case one should obtain only the real bounds similar to Lemma 7.4, then the results of Sections 8–10 can be applied immediately. This problem will be considered in a forthcoming paper.

The next theorem is the main result of the paper and Theorems B and C easily follow from it.

THEOREM A. *Let $f : X \hookleftarrow$ be a $C^3$ unimodal map of an interval to itself with a non-flat nonperiodic critical point $c$. Then there exists an interval $Z$ around the critical value $f(c)$ such that if $f^n(x) \in Z$ for $x \in X$ and $n > 0$, then $Sf^n(x) < 0$.*

Thus, if the orbit of some point passes nearby the critical value, the Schwarzian derivative becomes negative. However, one may prefer to work in a neighborhood of a critical point (not a critical value). In this case we do not get the negative Schwarzian derivative, but we have nice estimates for the cross-ratios.

THEOREM B. *Let $f : X \hookleftarrow$ be a $C^3$ unimodal map of an interval to itself with a non-flat nonperiodic critical point and suppose that the map $f$ does not have any neutral periodic points. Then there exists a constant $C_2 > 0$ such that if $M$ and $I$ are intervals, $I$ is a subinterval of $M$, $f^n|_M$ is monotone and $f^n(M)$ does not intersect the immediate basins of periodic attractors, then*

$$\begin{aligned} \mathbf{A}(f^n, M, I) &> \exp(-C_2 \, |f^n(M)|^2), \\ \mathbf{B}(f^n, M, I) &> \exp(-C_2 \, |f^n(M)|^2). \end{aligned}$$

Here a periodic attractor can be either a hyperbolic attracting periodic orbit or a neutral periodic orbit if its basin of attraction contains an open set. The *immediate basin* of a periodic attracting orbit is called a union of connected components of its basin which contain points of this orbit.

Notice the difference between Lemma 2.3 and this theorem: in Lemma 2.3 one has to estimate lengths of all intervals from the orbit of $M$ and in Theorem B one needs to know only the length of the last interval from the same orbit.

If the map does have a neutral repelling periodic point (i.e. a periodic point whose multiplier is $\pm 1$ and whose basin of attraction does not contain an open set) this theorem does not hold anymore. In this case the orbit of the interval $M$ can stay a very long time in the neighborhood of the neutral periodic orbit and the cross-ratio can became very small. (If the orbit of $M$ stays a long time in a neighborhood of some hyperbolic repelling orbit, we have some exponential expansion in this neighborhood and can control the sum of sizes of intervals from the orbit of $M$; then using Lemma 2.3 we can control the cross-ratios.) Fortunately, there is always a neighborhood of the critical point which does not contain any neutral periodic points. And of course, the most interesting dynamics is concentrated around the critical point.

There are many slightly different ways of generalizing the previous theorem to the case of maps with neutral periodic points. We will suggest a technical statement; however, it should cover all possible needs. But first we need to introduce some standard definitions.

We say that the point $x'$ is *symmetric* to the point $x$ if $f(x) = f(x')$. In this case we call the interval $[x, x']$ *symmetric* as well. A symmetric interval $I$ around a critical point of the map $f$ is called *nice* if its boundary points do not return into the interior of this interval under iterates of $f$. In the orbit of any periodic point one can take a point which is the nearest point in the orbit to the critical point. The interval between this point and its symmetric point will be nice. So there are nice intervals of arbitrarily small length if the critical point is not periodic.

Let $T \subset X$ be an interval and $f : X \hookleftarrow$ be a unimodal map. $R_T : U \to T$ denotes the first entry map to the interval $T$. The set $U$ consists of points whose iterates come to the interval $T$; i.e., $U = \{x \in X : \exists n > 0, \ f^n(x) \in T\}$. Now, if $x \in U$ and $n > 0$ is minimal such that $f^n(x) \in T$, then $R_T(x) = f^n(x)$. Notice that the set $U$ is not necessarily contained in the interval $T$. Sometimes we will want to consider only the points which are in the interval $T$ and in this case we will write $R_T|_T$ and the map $R_T|_T$ is called the first return map. So unless it is specifically mentioned otherwise, $R_T$ is defined on the set which can be larger than the interval $T$.

If the interval $T$ is nice, then the first entry map $R_T$ has some special properties. In this case the set $U$ is a union of intervals and if a connected component $J$ of the set $U$ does not contain the critical point of $f$, then $R_T : J \to T$ is a diffeomorphism of the interval $J$ onto the interval $T$. A connected component of the set $U$ will be called a *domain* of the first entry map $R_T$, or a *domain* of the nice interval $T$. If $J$ is a domain of $R_T$, the map $R_T : J \to T$ is called a *branch* of $R_T$. If a domain contains the critical point, it is called *central*.

THEOREM C. *Let $f$ be a $C^3$ unimodal map of an interval to itself with a non-flat critical point whose iterates do not converge to a periodic attractor. Then for any $0 < K < 1$ there is a nice interval $T$ around the critical point such that if*

- *$M$ is an interval and $f^n|_M$ is monotone,*

- *each interval from the orbit $\{M, f(M), \ldots, f^n(M)\}$ belongs to some domain of the first entry map $R_T$,*

*then*

$$\mathbf{A}(f^n, M, I) > K$$
$$\mathbf{B}(f^n, M, I) > K,$$

*where $I$ is any subinterval of $M$.*

The proofs of these theorems will occupy the rest of the paper.

## 4. The margins disjointness property

To be able to use Lemma 2.3 we need to bound the sum $\sum_{i=0}^{m-1} |f^i(M)|^2$ or the sum $\sum_{i=0}^{m-1} |f^i(M^-)| |f^i(M^+)|$. It is easy to do if the orbit of the interval $M$ is disjoint. Other useful estimates are formulated in the lemmas below.

In the previous version of this paper the results of this section played a key role in the proof that one can get nice cross-ratio bounds even if the map does not have negative Schwarzian derivative. The main theorem stated that the sum of squares of lengths of intervals from the orbit of some interval is small if the size of the last interval from this orbit is small. It appears that it is much easier to estimate the cross-ratios directly as is done in the present version of this paper. So the results of this section are used in the rest of the paper. However, I left them here because they can be useful, for example, in the multimodal case.

*Definition* 1. Assume a collection of oriented intervals $\{M_i\}$ and a collection of their subintervals $\{J_i\}$, $J_i \subset M_i$, $i = 1, \ldots, n$. Denote two components of the complement of $J_i$ in $M_i$ as $M_i^-$ and $M_i^+$ regarding the orientation. The collection $M_i \supset J_i$ has the *margins disjointness property* if and only if $M_i^- \cap M_j^- \neq \emptyset$ implies $M_i^+ \cap M_j^+ = \emptyset$ for $1 \leq i < j \leq n$.

Some forbidden and allowed configurations of intervals are shown in Figure 1.

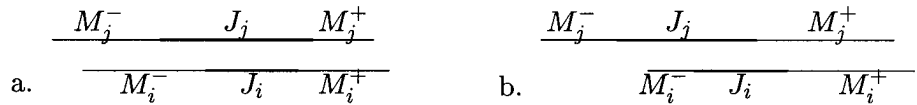

Figure 1. a. Forbidden configuration. b. Allowed configuration.

In the most important case we do not need to check the margin disjointness property for all pairs $(i, j)$ as is shown in the following obvious lemma:

LEMMA 4.1. *Let $f^n$ be strictly monotone on the interval $M$, $J \subset M$ and $M^-$ and $M^+$ are the components of the complement of $J$ in $M$. Then the collection $f^i(M) \supset f^i(J)$, $i = 0, \ldots, n$, has the margins disjointness property if and only if $f^i(M^-) \cap f^n(M^-) \neq \emptyset$ implies $f^i(M^+) \cap f^n(M^+) = \emptyset$ for $i = 0, \ldots, n-1$.*

LEMMA 4.2. *Let $X$ be an interval and $f : X \to X$ be a $C^0$ map. Let $M \supset J$ be intervals such that $f^n : M \to \tilde{M}$ is strictly monotone and the collection $\{f^i(M) \supset f^i(J), i = 0, \ldots, n\}$ has the margins disjointness property. Then*

$$\sum_{i=0}^n |f^i(M^-)|\,|f^i(M^+)| \leq 2|X| \max_{0 \leq i \leq n} |f^i(M)|.$$

◁ Let us consider rectangles $f^i(M^-) \times f^i(M^+)$, $i = 0, \ldots, n$ in the square $X \times X$. From the margins disjointness property it follows that these rectangles are pairwise disjoint. Moreover, they are contained in a narrow strip $S$ around the diagonal of the square

$$S = \{(x, y) \in X \times X : |x - y| \leq 2 \max_{0 \leq i \leq n} |f^i(M)|.$$

The sum we have to estimate is equal to the area of the union of the all rectangles and is bounded by the area of the strip.                                   ▷

The next lemma gives a simple way to check whether a collection of intervals has the margins disjointness property.

LEMMA 4.3.  *Let $f : X \to X$ be a $C^1$ map of the interval $X$. If $M \supset J$ are intervals such that $f^n$ is a diffeomorphism on $M$, $M \setminus J$ consists of two components $M^-$ and $M^+$ and $f^n$ maps $J$ on a critical point $c$ of the map $f$, $c \in f^n(J)$. Then the collection $\{f^i(M) \supset f^i(J), i = 0, \ldots, n\}$ has the margins disjointness property.*

◁ We will show that the condition $c \in f^n(J)$ implies the condition $f^i(M^-) \cap f^n(M^-) \neq \emptyset \Rightarrow f^i(M^+) \cap f^n(M^+) = \emptyset$. Indeed, assume that this is not true and there is the integer $i$, $i < n$, such that $f^i(M^-)$ intersects $f^n(M^-)$ and $f^i(M^+)$ intersects $f^n(M^+)$. Then $f^i(M)$ covers the whole interval $f^n(J)$. The critical point $c$ is in $f^n(J)$, so the map $f^n|_M$ is not a diffeomorphism.   ▷

## 5. Consequences of absence of wandering intervals

The interval $J$ is called a *wandering* interval of the map $f$ if it satisfies two conditions:

- The intervals of the forward orbit $\{f^i(J), i = 0, 1, \ldots\}$ are pairwise disjoint;
- The images $f^i(J)$ do not converge to a periodic attractor with $i \to \infty$.

Fortunately, in our case we will not have wandering intervals due to the following well-known theorem (see [3dMvS]):

THEOREM.  *If $f$ is a $C^2$ map with non-flat critical points, then $f$ has no wandering intervals.*

However, we will use not this theorem itself but its simple corollaries.

LEMMA 5.1.  *Let $f$ be a $C^2$ map with non-flat critical points and $J$ be an interval. Then either there is $n$ such that $f^n|_J$ is not monotone or the iterates of all points of the interval $J$ converge to some periodic orbits.*

◁ Assume the contrary, i.e. assume that $f^i|_J$ is monotone for any $i > 0$ and that the iterates of $J$ do not converge to a periodic attractor. Consider the set $U = \bigcup_{i=0}^{\infty} f^i(J)$ and take the connected component $T$ of this set that contains $J$. The set $U$ is forward invariant and $f^i|_T$ is monotone. If $f^i(T) \cap f^j(T) \neq \emptyset$ for $i < j$, then $f^i(T) \supseteq f^j(T)$. Since $f^{j-i} : f^i(T) \to f^j(T)$ is monotone all points of $T$ will converge to some periodic orbit. The other possibility is that the orbit of $T$ is disjoint and therefore $T$ is a wandering interval. We arrived at a contradiction in both cases. ▷

LEMMA 5.2. *Let $f$ be a $C^2$ map with non-flat critical points. Then there is a function $\tau_1$ such that $\lim_{\varepsilon \to 0} \tau_1(\varepsilon) = 0$ and such that if $V$ is an interval, $f^n|_V$ is a diffeomorphism, and $f^n(V)$ is disjoint from the immediate basins of periodic attractors, then*

$$\max_{0 \leq i \leq n} |f^i(V)| < \tau_1(|f^n(V)|).$$

◁ Suppose that such a function $\tau_1$ does not exist, i.e. there is a constant $\varepsilon > 0$ such that for any $\delta > 0$ there is an interval $V$ of length greater than $\varepsilon$ and such that $|f^n(V)| < \delta$ for some $n$. Moreover, the map $f^n|_V$ is a diffeomorphism and the interval $f^n(V)$ does not intersect the immediate basin of attraction. Take a sequence of $\delta_i$ tending to 0 and sequences of corresponding intervals $V_i$ and corresponding iterates $n_i$. Extract a convergent subsequence $V_{i_j}$ and denote its limit as $V_0$. The interval $V_0$ cannot be degenerate because its length is greater than or equal to $\varepsilon$. The sequence $n_{i_j}$ tends to infinity; otherwise we could take a bounded subsequence and we would have that $f^{n_0}(V_0)$ is a point and this is impossible. The maps $f^{n_i}|_{V_0}$ are diffeomorphisms and the interval $f^{n_i}(V_0)$ does not intersect the immediate basin of attraction. This contradicts the previous lemma. ▷

## 6. High, low and center returns

As mentioned before, the most interesting dynamics of a unimodal map is concentrated in the neighborhood of its critical point, so it is natural to consider a first return map to some neighborhood of the critical point. Then we can observe dynamics of the map as under "a microscope". The first return map to a nice interval has particularly nice properties: the boundary of domains of the first return map is mapped to the boundary of the nice interval.

Let $f$ be a unimodal map, $T$ be a nice interval, $R_T$ be the first entry map to $T$ and $J$ be its central domain. We will need to distinguish between different types of the first entry maps depending on the position of the image of the critical point. First, if the image under $R_T$ of the central domain covers the critical point, then $R_T$ is called a *high* return; otherwise it is a *low* return. And if $R_T(c)$ is in the central domain $J$, where $c$ is a critical point of $f$, then $R_T$ is a *central* return (otherwise it is noncentral).

It is possible that the first return map $R_T|_T$ will have just one unimodal branch defined on the whole interval $T$. If this happens, such a map is called *renormalizable* and $T$ is called a *restrictive* interval. If for a map there exists a sequence of restrictive intervals, then this map is called *infinitely renormalizable*.

## 7. Extensions of branches

In this section we will construct some space around nice intervals and prove that the range of some branches of the first entry map to this interval can be extended to this space.

Suppose that $g : X \hookleftarrow$ is a $C^1$ map and suppose that $g|_V : V \to J$ is a diffeomorphism of the interval $V$ onto the interval $J$. If there is a larger interval $V' \supset V$ such that $g|_{V'}$ is a diffeomorphism, then we will say that the range of the map $g|_V$ can be *extended* to the interval $g(V')$.

LEMMA 7.1. *Let $f$ be a unimodal map, $T$ be a nice interval, $J$ be its central domain and $V$ be a domain of the first entry map to $J$ which is disjoint from $J$, i.e. $V \cap J = \emptyset$. Then the range of the map $R_J : V \to J$ can be extended to $T$.*

◁ Let $I \supset V$ be a maximal interval of monotonicity of $f^n$ where $f^n = R_J|_V$. By some iterations of $f$ the boundary points of $I$ are mapped on the critical point while the image of the interval $V$ stays outside of $J$. So a boundary point of $J$ belongs to some iterate of the interval $I$. $J$ is the central domain of $T$; thus the boundary points of $J$ will never return inside $T$. This implies that the interval $f^n(I)$ covers the whole interval $T$. ▷

The next lemma deals with the renormalizable and almost renormalizable cases and uses the method of the smallest interval, which was used by Martens to prove a similar statement for renormalizable maps.

LEMMA 7.2. *Let $f$ be a $C^3$ unimodal map with a non-flat recurrent critical point $c$. There exist constants $0 < \tau_2 < 1$ and $\tau_3 > 0$ such that if $T$ is any sufficiently small nice interval around the critical point $c$ and its central domain $J$ is sufficiently big, i.e. $\frac{|J|}{|T|} > \tau_2$, then there is an interval $W$ containing a $\tau_3$-scaled neighborhood of the interval $T$ such that*

- *if $c \in R_T(J)$ (i.e. $R_T$ is a high return), then the range of any branch $R_T : V \to T$ can be extended to $W$ provided that the domain $V$ is not contained in $T$;*

- if $c \notin R_T(J)$ (i.e. $R_T$ is a low return) and the map $f$ is not renormalizable, then the range of the branch $R_T : J_1 \to T$ of the first entry map to the interval $T$ can be extended to $W$, where $J_1$ is a domain of $R_T$ containing the critical value $f(c)$.

◁ First we will construct some space around the interval $T$.

Let the central branch of the first entry map have the form $f^k : J \to T$. Consider the orbit $\{f^i(J), i = 1, \ldots, k\}$ of the interval $J$. Since $f^k|_J$ is the first entry map, the orbit of the interval $J$ is disjoint. Take an interval of minimal length in this orbit. Denote this interval as $U$ and denote the interval which is a 1-scaled neighborhood of $U$ as $M$. If the interval $T$ is sufficiently small, then the interval $M$ lies in the domain of definition of the map $f$. The pullback of the interval $M$ along the orbit of $J$ will give us a required space around the interval $T$.

First, the interval $M$ does not contain any other intervals from the orbit $\{f^i(J), i = 0, \ldots, k\}$ because of the minimality of $U$. (However, $M$ can have non empty intersections with two intervals from the orbit different from $U$.)

The range of the map $f^{k_1-1} : f(J) \to U$ can be diffeomorphically extended to $M$. Indeed, denote the maximal interval of monotonicity of $f^{k_1-1}$ containing the interval $f(J)$ as $I$. The boundary points $a_-$ and $a_+$ of the interval $I$ are mapped on the critical point $c$ by some iterates of $f$, $f^{i_\pm}(a_\pm) = c$, $i_-, i_+ < k_1$. The interval $f^{i_\pm+1}(J)$ is outside of $T$, hence on both ends of the interval $f^{k_1}(I)$ there are intervals from the orbit of $J$ different from $U$. This implies that $M$ belongs to $f^{k_1-1}(I)$.

The pullback $\{M_i, i = 0, \ldots, k_1\}$ of the interval $M$ along the orbit $\{f^i(J), i = 0, \ldots, k_1\}$ has intersection multiplicity bounded by 4 (this means that any point of the interval $X$ belongs to at most four intervals $M_i$, $i = 0, \ldots, k_1$). Indeed, suppose that this is not true and there is a point $b$ which belongs to five intervals from the pullback of $M$. Then there are three intervals from those five, which we denote $\{M_{i_1}, M_{i_2}, M_{i_3}\}$, such that $f^{i_2}(J) \subset M_{i_1}$, $f^{i_3}(J) \cap M_{i_1} \neq \emptyset$ and the interval $f^{i_2}(J)$ is situated between the intervals $f^{i_1}(J)$ and $f^{i_3}(J)$ (see Fig. 2). We know that the intervals $J_1, f(J_1), \ldots, f^{k-1}(J_1) = T$ are disjoint because $f^{k-1} : J_1 \to T$ is a branch of the first entry map; hence $f^{i_2-1}(J_1) \subset M_{i_1}$. If $k_1 - i_1 \leq k - i_2$, then $k \geq i_2 + k_1 - i_1$ and $f^{i_2+k_1-i_1}(J) \subset M = M_{k_1} = f^{k_1-i_1}(M_{i_1})$, this is a contradiction to the choice of the interval $M$. If $k_1 - i_1 > k - i_2$, then $k_1 > i_1 + k - i_2$ and applying the map $f^{k-i_2}$ to the intervals $f^{i_2-1}(J_1) \subset M_{i_1}$ we obtain the inclusion $T \subset M_{i_1+k-i_2}$. So the critical point $c$ is contained in the interval $M_{i_1+k-i_2}$, this contradicts the fact that the map $f^{k_1-1} : M_1 \to M$ is a diffeomorphism.

We can apply Lemma 2.3 and obtain some definite space around the interval $f(J)$. The critical point $c$ of the map $f$ is nonflat, so we can pull back the space to the interval $J$. Indeed, near the critical point the map $f$ has the form $f(x) = |\phi(x)|^\alpha + f(c)$. Now if we have two points $a$ and $b$ which belong to the domain of definition of $\phi$ and such that $b \in [a, c]$, we obtain the inequality

$$\frac{|b-c|}{|a-b|} < C_3 \left(\frac{|f(b) - f(c)|}{|f(a) - f(b)|}\right)^{1/\alpha}$$

where $C_3 > 0$ depends only on the diffeomorphism $\phi$. So the maximal interval $W$ around $J$ which maps onto the interval $M$ by $f^{k_1}$ is a $\tau_4$-scaled neighborhood of $J$ where $\tau_4$ is some universal constant.

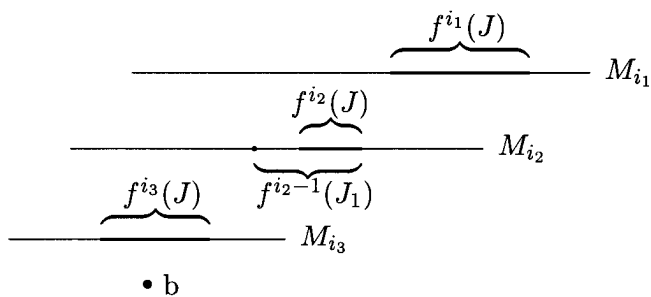

Figure 2. Three overlapping intervals.

Instead of counting the intersection multiplicity of the intervals $\{M_i\}$ we could observe that the collection of intervals $\{M_i \supset f^i(J), i = 1, \ldots, k_1\}$ has the margins disjointness property (this is trivial) and then Lemmas 4.2 and 2.3 immediately imply that the interval $J_1$ has some definite space inside the interval $M_1$.

If we choose the constant $\tau_2$ to be sufficiently close to 1, then the interval $W$ will cover the interval $T$ and thus $W$ will give a definite space around the interval $T$.

Suppose that $R_T$ is a high return, so that $c \in R_T(J)$, and let $V$ be a domain of $R_T$ not containing the critical point, i.e. $V \neq J$. We want to show that the range of the map $f^{k_2} : V \to T$ can be extended to the interval $W$, where $f^{k_2} = R_T|_V$. Indeed, arguing as before we can conclude that the range of the map $f^{k_2} : V \to T$ can be extended to an interval which contains the interval $T$ and two intervals of the form $f^j(J)$ with $0 < j < k$ on either side of $T$. If these intervals were contained in $W$, then the interval $M$ would contain the interval $f^{j_1}(J)$, where $j_1 \equiv k_1 + j \pmod{k}$.

Now let us prove that the range of the map $f^{k-1} : J_1 \to T$ can be diffeomorphically extended to the interval $W$ even if $R_T$ is a low return. If this is not the case, then there is an interval $\hat{W}$ which is a pullback of $W$ to the critical value $f(c)$ along the orbit of the interval $J$ such that $f(c) \in \hat{W}$, $f^{j-1}(J_1) \subset \hat{W}$, $f^{k-j}(\hat{W}) = W$ and $f^{k-j}|_{\hat{W}}$ is a diffeomorphism, where $1 < j < k$ (see Fig. 3). The interval $f(W)$ cannot contain the interval $f^{j-1}(J_1)$. Indeed, if

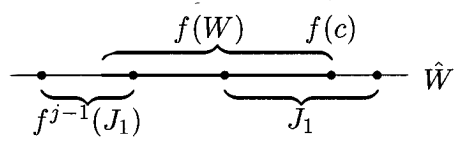

Figure 3. The interval $f(W)$ cannot contain the interval $f^{j-1}(J_1)$, and therefore it is in $\hat{W}$.

$k_1 - 1 < k - j$, then $f^{j+k_1-1}(J) \subset f^{j+k_1-2}(J_1) \subset f^{k_1}(W) \subset M$; if $k_1 - 1 \geq k - j$, then $c \in T = f^{k-1}(J_1) \subset f^{k-j+1}(W) \subset M_{k-j+1}$. Both cases are impossible and $f^{j-1}(J_1) \not\subset f(W)$. Thus, $f(W) \subset \hat{W}$. This implies that the map $f$ is renormalizable and $f^{k-j+1}(W) \subset W$. This contradicts the assumptions on the map $f$. ▷

LEMMA 7.3. *Let $f$ be a $C^3$ unimodal map with a non-flat recurrent critical point $c$. There are constants $\tau_5 < 1$ and $C_4 > 0$ such that if $T$ is a nice interval, $|T| < C_4$, the first entry map $R_T$ is a non-central low return and $J$ is a central domain of $R_T$, then*
$$\frac{|J|}{|T|} < \tau_5.$$

◁ Suppose that $\frac{|J|}{|T|} > \tau_2$ (otherwise we have nothing to do). According to the previous lemma there are the interval $W$ which is a $\tau_3$–scaled neighborhood of the interval $T$ and an interval around the critical value $f(c)$ which is diffeomorphically mapped onto the interval $W$. Denote this latter interval as $U_1$, so that $f^{k-1}(U_1) = W$, and let $U$ be the full preimage of $U_1$ under the map $f$. Suppose also that $f^k(c) < c$ (this is not a restriction) and let $R$ be equal to $T \setminus f^k(J)$, $L$ be a component of $W \setminus T$ such that the interval $f^k(J)$ is situated between the intervals $L$ and $R$, and let $T' = T \cup L$ (see Fig. 4).

If $r$ is a pullback of the interval $R$ under the map $f^{k-1}$, then the orbit $\{f^i(r), i = 0, \ldots, k-1\}$ is disjoint because $r \subset J_1$ and the orbit of the interval $J_1$ is disjoint. Hence the sum $\sum_{i=0}^{k-1} |f^i(l)||f^i(r)|$ is small if the interval $T$ is small (here $l$ is a pullback of $L$). Applying Lemma 2.3 we obtain the inequality
$$\mathbf{a}(l \cup J_1, f(J)) < C_5 \, \mathbf{a}(T', f^k(J)),$$
where the constant $C_5$ is close to 1 if the interval $T$ is small.

If the interval $f^k(J)$ was very small compared with the interval $T$ (and this is inevitable if the return is noncentral low and $\frac{|J|}{|T|}$ is close to one), then the ratio $\frac{|f(J)|}{|l \cup f(J)|} < \mathbf{a}(l \cup J_1, f(j))$ would be very small and as a consequence the ratio $\frac{|J|}{|U|}$ would be very small as well. Therefore the interval $U$ would be much larger than the intervals $J$, $T$ and $W$. In this case we would have the

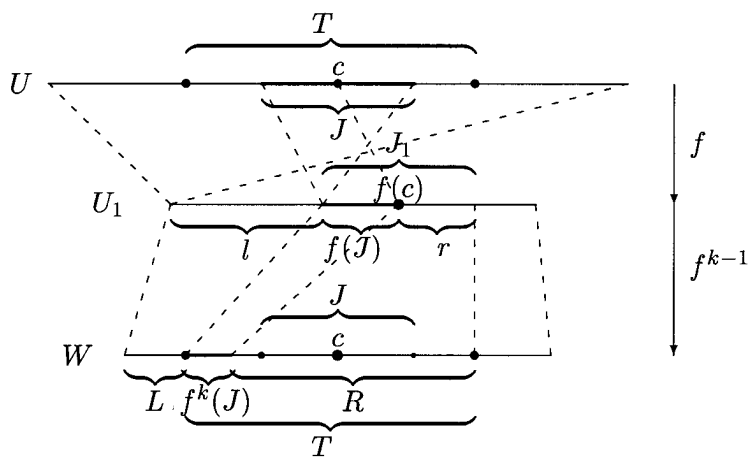

Figure 4. Iterates of the interval $J$.

unimodal map $f^k : U \to U$ such that $c \notin f^k(U)$. This would imply that the iterates of the critical point $c$ converge to some periodic attractor. This is a contradiction to the assumption that $c$ is recurrent, and hence the ratio $\frac{|J|}{|T|}$ cannot be close to one. $\triangleright$

Here we summarize the previous three lemmas:

LEMMA 7.4. *Let $f$ be a $C^3$ unimodal map with a non-flat nonperiodic critical point. There is a constant $\tau_6$ and a sequence $\{T_i, i = 1, \ldots\}$ of nice intervals whose sizes shrink to $0$ such that the range of any branch $R_{T_i} : V \to T_i$ of the first entry map can be extended to an interval which contains a $\tau_6$-scaled neighborhood of $T_i$ provided that the domain $V$ is disjoint from $T_i$.*

$\triangleleft$ If the critical point is not recurrent, then the statement is obvious. Indeed, if $c$ is not recurrent, then there exists an interval $W$ around $c$ which does not contain any other points of the forward orbit of $c$. If $T$ is any nice interval contained in $W$, then the range of any branch of $R_T$ can be extended to $W$ (the proof is the same as the proof of Lemma 7.1) and the lemma follows. So we will assume that $c$ is recurrent.

Let us consider several cases. First, suppose that the map $f$ is infinitely renormalizable. Then a sequence $\{T_i\}$ is just a sequence of restrictive intervals (possibly we will have to drop the beginning of the sequence of the restrictive intervals in order that these intervals become very small and Lemma 7.2 starts to work because the first entry map to a restrictive interval is always a high return).

Now let the map $f$ be only finitely renormalizable. Take any small nice interval $T'_1$ and consider a sequence of intervals $T'_1, T'_2, \ldots$ such that the interval $T'_{i+1}$ is a central domain of the interval $T'_i$. If the interval $T'_1$ is taken to be sufficiently small, the lengths of the intervals $T'_i$ will tend to zero. If in the sequence $\{R_{T'_i}\}$ there are infinitely many high returns $\{R_{T'_{i_j}}\}$, then

$$T_j = \begin{cases} T'_{i_j}, & \text{if } \frac{|T'_{i_j+1}|}{|T'_{i_j}|} > \tau_2 \\ T'_{i_j+1}, & \text{if } \frac{|T'_{i_j+1}|}{|T'_{i_j}|} \leq \tau_2 \end{cases}.$$

If there are only low returns, then there exist infinitely many noncentral low returns $R_{T'_{i_j}}$ (otherwise the critical point $c$ would be nonrecurrent). In this case we put

$$T_j = T'_{i_j+1}. \qquad \triangleright$$

## 8. Derivative estimate

We cannot use the Koebe principle and estimate the distortion of some iterate of the map $f$ restricted to some interval $J$ if we do not know a bound of the sum of squares of lengths of intervals from the orbit of the interval $T$ whose image is definitely larger than the image of $J$. However, it appears that we can estimate the derivative from below if we can bound the sum of lengths of intervals from the orbit of $J$.

LEMMA 8.1. *Let $f : X \hookleftarrow$ be a $C^3$ map with non-flat critical points and let $J \subset T$ be intervals such that $f^n|_T$ is monotone, the interval $f^n(T)$ contains a $\delta$-scaled neighborhood of the interval $f^n(J)$ and the orbit $\{f^i(J), i = 0, \ldots, n-1\}$ is disjoint. Then there exists a constant $C_6 > 0$ depending only on the map $f$ such that*

$$|Df^n(x)| > C_6 \frac{\delta}{1+\delta} \frac{|f^n(J)|}{|J|}$$

*where $x \in J$.*

◁ Let the point $x$ cut the interval $J$ onto two intervals $J^-$ and $J^+$. Obviously, one of the following two inequalities must hold: either $\frac{|f^n(J^-)|}{|J^-|} > \frac{|f^n(J)|}{|J|}$ or $\frac{|f^n(J^+)|}{|J^+|} > \frac{|f^n(J)|}{|J|}$ (see Fig 5). Suppose that the second inequality (for $J^+$) holds.

Let $J'$ be an infinitesimal interval around the point $x$ and let $T' = T^- \cup J$. If we apply Lemma 2.3 to the intervals $J' \subset T'$ and the map $f^n$, then we obtain the following inequality:

$$\begin{aligned} \mathbf{A}(f^n, T', J') &> \exp\left(C_1 \max_{0 \leq i \leq n} |f^i(T^- \cup J^-)| \sum_{i=0}^{n} |f^i(J^+)|\right) \\ &> \exp\left(C_1 \max_{0 \leq i \leq n} |f^i(T)| \sum_{i=0}^{n} |f^i(J)|\right). \end{aligned}$$

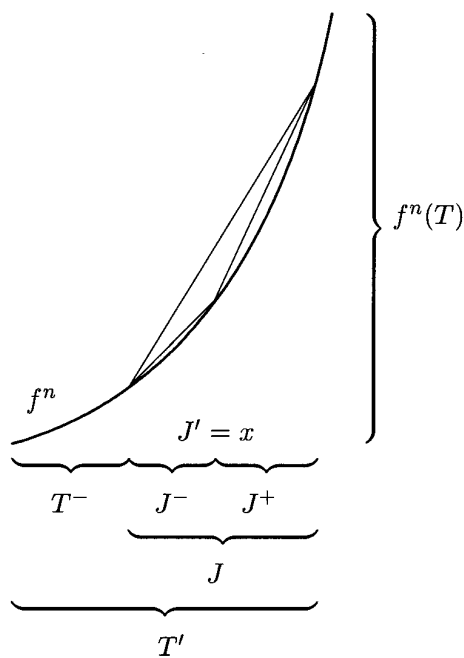

Figure 5. One of the slopes is greater than the average slope.

Notice that $\max_{0 \leq i \leq n} |f^i(T)| \sum_{i=0}^{n} |f^i(J)| < 2|X|$ because the orbit of $J$ is disjoint and put $C_6 = \exp(-2|X|C_1)$.

The ratio $\frac{|f^n(J')|}{|J'|}$ is just the derivative $|Df^n(x)|$. So, rearranging the terms in the previous inequality we get:

$$
\begin{aligned}
|Df^n(x)| &> C_6 \frac{|T'|}{|J^+|(|T^-|+|J^-|)} \frac{|f^n(J^+)|(|f^n(T^-)|+|f^n(J^-)|)}{|f^n(T')|} \\
&> C_6 \frac{|f^n(J^+)|}{|J^+|} \frac{|f^n(T^-)|+|f^n(J^-)|}{|f^n(T')|} \\
&> C_6 \frac{\delta}{1+\delta} \frac{|f^n(J)|}{|J|}. \qquad \triangleright
\end{aligned}
$$

## 9. The Schwarzian derivative of the first entry map

*Proof of Theorem* A. First, let us consider the case when the trajectory of the critical point does not converge to a periodic attractor.

Let $T$ be a nice interval around $c$ from the sequence given by Lemma 7.4 and let $T$ be so small that $T$ is disjoint from the immediate basins of attractors. Let $f^n : V \to T$ be a branch of the first entry map to $T$ and $V \not\subset T$. As we know, the map $f^n : V \to T$ is a diffeomorphism and its range can be diffeomorphically extended to $W$ where the interval $W$ contains a $\tau_6$-scaled neighborhood of $T$.

Due to Lemma 8.1 we can estimate the derivative of $f^{n-i} : f^i(V) \to T$ by the ratio of intervals: $|Df^{n-i}(x)| > C_7 \frac{|T|}{|f^i(V)|}$, $x \in f^i(V)$, $i = 0, \ldots, n$, where $C_7 = C_6 \frac{\tau_6}{1+\tau_6}$.

The map $f$ has a nonflat critical point and nearby the critical point has the form $f(x) = \pm|\phi(x)|^\alpha + f(c)$, where $\phi$ is some local $C^3$ diffeomorphism with $\phi(c) = 0$, $\alpha \geq 2$. The Schwarzian derivative of the function $x^\alpha$ is equal to $S(x^\alpha) = \frac{1-\alpha^2}{2x^2}$ and since $\phi$ is a diffeomorphism its Schwarzian derivative is bounded by some constant, $|S\phi(x)| < C_8$. Thus, applying the Schwarzian derivative to the composition of the functions $\phi$ and $x^\alpha$ we obtain

$$Sf(x) = \frac{1-\alpha^2}{2\phi(x)^2}(D\phi(x))^2 + S\phi(x).$$

Hence, if $T$ is sufficiently small, there exists a constant $C_9 > 0$ such that $Sf(x) < -\frac{C_9}{(x-c)^2}$ for $x \in T$. Outside of the interval $T$ the map $f$ has no critical points, therefore the Schwarzian derivative of $f$ is bounded there by some constant $C_{10} > 0$, i.e. $|Sf(x)| < C_{10}$ for $x \notin T$. Since near the critical point the Schwarzian derivative of $f$ is negative, $Sf(x) < C_{10}$ for all $x \in X$.

Now let us estimate the Schwarzian derivative of the map $f^{n+1} : V \to f(T)$.

$$\begin{aligned} S(f^{n+1})(x) &= Sf(f^n(x))\,|Df^n(x)|^2 + \sum_{i=0}^{n-1} Sf(f^i(x))\,|Df^i(x)|^2 \\ &= |Df^n(x)|^2 \left( Sf(f^n(x)) + \sum_{i=0}^{n-1} Sf(f^i(x))\,|Df^{n-i}(f^i(x))|^{-2} \right) \\ &\leq \left(\frac{|Df^n(x)|}{|T|}\right)^2 \left( -C_9 \left(\frac{|T|}{f^n(x)-c}\right)^2 + C_7 C_{10} \sum_{i=0}^{n-1} |f^i(V)|^2 \right). \end{aligned}$$

Note that $\frac{|T|}{|f^n(x)-c|}$ is always greater than 1 because $f^n(x) \in T$. The intervals from the orbit of $V$ are disjoint; thus

$$\sum_{i=0}^{n-1} |f^i(V)|^2 < |X| \max_{0 \leq i < n} |f^i(V)| < |X|\tau_1(|T|)$$

(see Lemma 5.2). As a result we have that if the nice interval $T$ is small enough, then the first entry map to $f(T)$ has negative Schwarzian derivative. If $f^m(y) \in f(T)$ for some $y \in X$, $n > 0$, then $f^m$ can be decomposed as $R_{f(T)}^{m'}$. Each branch of $R_{f(T)}$ has negative Schwarzian derivative; thus $S(f^m)(y) < 0$.

Now consider the case when the trajectory of $c$ is attracted to some periodic orbit. (In fact, this is not really an interesting case because the dynamics of such maps is very well understood.)

The estimate for $S(f^{n+1})(x)$ given above is still valid, but we cannot use Lemma 5.2 any more. However the sum $\sum_{i=0}^{n-1} |f^i(V)|^2$ is still bounded by $|X|^2$.

Notice that if the interval $T$ is disjoint from $f^k(T)$ for all $k > 0$, then the first entry map to any subinterval $T'$ of $T$ is a restriction of the first entry map to $T$, i.e., $R_{T'} = R_T|_{U'}$. So if $T'$ is very small and $f^n(x) \in T'$, then the ratio $\frac{|T|}{|f^n(x)-c|}$ is very large; hence, $Sf^n(x) < 0$ and $SR_{f(T')} < 0$. □

## 10. Cross-ratio estimate

*Proof of Theorem* B. Take a nice interval $T$ around the critical point $c$ such that the first entry map to $f(T)$ has negative Schwarzian derivative. Denote by $J$ a symmetric interval around $c$ which is one-half the size of $T$, $|J| = |T|/2$ and denote two components of the complement of the interval $J$ in the interval $T$ as $T^-$ and $T^+$.

Divide the orbit of $M$ into 3 parts. Let $n_1$ be the maximal integer satisfying the following property: $f^{n_1-1}(M) \subset T$. Then the first (maybe empty) part of the orbit is $\{f^i(M), i = 0, \ldots, n_1\}$. The map $f^{n_1}|_M$ has negative Schwarzian derivative, hence it can only increase the cross ratio.

Let $n_2$ be the minimal integer satisfying the following property: $T^-$ or $T^+$ belongs to the interval $f^{n_2}(M)$. Then the second part of the orbit consists of the intervals $\{f^i(M), i = n_1 + 1, \ldots, n_2\}$ and the third part is the rest. If such a number $n_2$ does not exist, then we put $n_2 = n$.

It is easy to see that all but the last interval in the second part lie outside of the interval $J$. Due to the theorem of Mañé ([Man]) there are constants $\tau_7 < 1$ and $C_{11} > 0$ such that $|f^i(M)| < C_{11} \tau_7^{n_2-i} |f^{n_2}(M)|$. Combining this estimate and the claim above we obtain the bound on the sum of squares of lengths of intervals in the second part of the orbit: $\sum_{i=n_1+1}^{n_2} |f^i(M)|^2 < C_{12} |f^{n_2}(M)|^2$.

Since there are no wandering intervals, there exists an integer $N$ such that $f^{N+1}|_{T^\pm}$ is not monotone. Let $N$ be minimal with this property; then the number of intervals in the third part is bounded by $N$. The following claim is obvious:

*Claim.* There is a constant $C_{13}$ such that if $V$ is an interval, $f^k|_V$ is a diffeomorphism, where $k \leq N$, $f^i(V)$ does not belong to the interval $T$ for $i = 0, \ldots, k$, then $|V| < C_{13} |f^k(V)|$.

Thus, $\sum_{i=n_2}^{n} |f^i(M)|^2 < NC_{13}|f^n(M)|^2$. Applying Lemma 2.3 we complete the proof. □

*Proof of Theorem* C. This theorem is a simple consequence of Theorem A. Indeed, let the nice interval $T$ be so small that the first entry map $R_{f(T)}$ has negative Schwarzian derivative. Let us decompose the map $f^n$ as $f^k \circ R_{f(T)}^j$ and let $k$ be the smallest positive integer with such a property. $R_{f(T)}^j$ has negative Schwarzian derivative, so it does not decrease the cross-ratios. On the other hand,

$$B(f^k, f^{n-k}(M), f^{n-k}(I)) > \exp(C_1 \sum_{m=n-k}^{n} |f^m(M)|^2).$$

The intervals $f^{n-k}(M), f^{n-k+1}(M), \ldots, f^n(M)$ are disjoint because they are contained in the domains of the first entry map $R_T$, so the sum can be bounded:

$$\sum_{m=n-k}^{n} |f^m(M)|^2 < |X| \max_U |U|,$$

where $U$ runs over all domains of the first entry map $R_T$. The last quantity can be made arbitrarily small. □

MATHEMATICS INSTITUTE, UNIVERSITY OF WARWICK, COVENTRY, UK
*E-mail address*: oleg@maths.warwick.ac.uk